
\input amstex
\magnification 1200
\documentstyle{amsppt}
\topmatter
\title Rationality of the Moduli Variety of Curves of Genus 3
\endtitle
\author P.Katsylo
\endauthor
\address Independent University of Moscow
\endaddress
\date June 25, 1994
\enddate
\thanks Research supported by Max-Planck-Institute fur Mathematik
  and Grant N MQZ000 of the international Science Foundation
\endthanks
\abstract We prove rationality of the moduli variety of curves of
  genus 3.
\endabstract
\endtopmatter
\document
\define\C{\Bbb C}
\define\p{\Bbb P}
\par
\S\bf 0. \rm Let $g \ge 2$ be a natural number.
Consider the moduli variety $M_g$
of curves of genus $g$. Recall that $M_g$ is an irreducible
quasiprojective variety, $\dim M_g = 3g - 3$ \cite{1,2}.
For $g \ge 23 \ \ M_g$ is not unirational \cite{3}.
If $g \le 13$, then $M_g$ is unirational \cite{4-6}.
For $g = 2, 4 ,5 ,6 \ \ \ M_g$ is rational \cite{7-10}.
The aim of this paper is to prove the following result.

\proclaim{Theorem 0.1} $M_3$ is rational.
\endproclaim

\par
The group $SL_3$ acts canonically in the space $S^4\C^{3\ast}$.
As is known,
$$\C(M_3) \approx \C(P(S^4\C^{3\ast}))^{SL_3}. \tag 0.1
$$
\par
For $d \ge 0$ denote by $V(2d)$ the space of forms of degree $2d$
in the variables $z_1, z_2$. The group $PSL_2$ acts canonically
in $V(2d)$. For
$\lambda = (\lambda_0, \lambda_2, \lambda_4, \lambda_6) \in \C^{\ast 4}$
consider the homogeneous (of degree 2) $PSL_2$-morphism
$$\gathered
    \delta_{\lambda} : V(8) \dotplus V(0) \dotplus V(4)
     \rightarrow V(4), \\
      f_8 + f_0 + f_4 \mapsto \lambda_6
    \psi_6(f_8, f_8) + \lambda_4 \psi_4(f_8, f_4)
      + \lambda_2 \psi_2(f_4, f_4) + \lambda_0 f_4 f_0.
  \endgathered
$$
Here $\psi_i$ is \it i\rm th \ transvectant.
Consider $\delta^{-1}_{\lambda}(0)$. It is obvious that
$1 \in \delta^{-1}_{\lambda}(0)$ and the tangent space to
$\delta^{-1}_{\lambda}(0)$ at the point $1$ coinsides with
$V(8) \dotplus V(0)$. It follows that $1$ is a regular point
of the subvariety $\delta^{-1}_{\lambda}(0)$. Therefore, a unique
($10$-dimensional) irreducible component $U_{\lambda}$ of the
subvariety $\delta^{-1}_{\lambda}(0)$ contains $1$. We have the
following isomorphism of the fields
$$\C(P(S^4\C^{3\ast}))^{SL_3} \approx
  \C(U_{(-7/36,11/54,1/840,-6/1225)})^{PSL_2 \times \C^{\ast}} \tag 0.2
$$
(see \cite{11}).
\proclaim{Theorem 0.2} For all $\lambda \in \C^{\ast 4}$ the field
  $\C(U_{\lambda})^{PSL_2 \times \C^{\ast}} \approx
  \C(PU_{\lambda})^{PSL_2}$ is rational.
\endproclaim
\par
Theorem 0.1 is a consequence of (0.1), (0.2), and Theorem 0.2. We prove
Theorem 0.2 in \S\S 1-6.
\par
This paper is organized as follows. In \S 1 we reduce
Theorem 0.2 to its partial case
$\lambda = (1, 6\epsilon, 1, 6)$, where $\epsilon \neq 0$.
Then we fix a basis $e_1, \dots , e_9, a_0, \dots , a_5$
in the space $V(8) \dotplus V(0) \dotplus V(4)$ and
write down the mapping $\delta_{\lambda}$ in coordinates.
In \S 2 we recall some facts about $(G,G')$-sections.
In \S 3 we construct $(PSL_2, N(H))$-section
$PX_{\lambda}^0$ of the variety $PU_{\lambda}$.
We have the isomorphisms of the fields
$$\C(PU_{\lambda})^{PSL_2} \approx
  \C(PX_{\lambda}^0)^{N(H)} \approx
  \C(PX_{\lambda})^{N(H)},
$$
where $X_{\lambda} = \overline{X_{\lambda}^0}$.
In \S 4 we construct $6$-dimensional variety
$Y_{\lambda}$ and regular action $N(H) : Y_{\lambda}$
such that
$$\C(PX_{\lambda})^{N(H)} \approx
  \C(Y_{\lambda})^{N(H)}
$$
and the subgroup $H \subset N(H)$ acts on
$Y_{\lambda}$ trivially.
In \S\S 5,6 we prove rationality of the field
$\C(Y_{\lambda})^{N(H)}$.
\par
The author is grateful to \`E.B.Vinberg,
V.A.Iskovskikh and S.L.Tregub for useful discussions.
\par
\S \bf 1. \rm Note that it is sufficient to prove theorem 0.2 for
$\lambda = (1, 6\epsilon, 1, 6)$, where $\epsilon \neq 0$. Indeed,
suppose $6\mu^2_8 = \lambda_6,\ \mu_4\mu_8 = \lambda_4,
\ 6\epsilon\mu^2_4 = \lambda_2,\ \mu_0\mu_4 = \lambda_0$; then
$$PU_{(\lambda_0, \lambda_2, \lambda_4, \lambda_6)} \rightarrow
  PU_{(1, 6\epsilon, 1, 6)}, \quad
  \overline{\mu_8 f_8 + \mu_0 f_0 + \mu_4 f_4} \tag 1.1
$$
is $PSL_2$-isomorphism. From (1.1) it follows that
$$\C(PU_{(\lambda_0, \lambda_2, \lambda_4, \lambda_6)})^{PSL_2}
  \approx \C(PU_{(1, 6\epsilon, 1, 6)})^{PSL_2}.
$$
\par
Let us prove Theorem 0.2 for $\lambda = (1, 6\epsilon, 1,6)$,
where $\epsilon \neq 0$.
\par
Fix the following basis in the space
$V(8) \dotplus V(0) \dotplus V(4)$ :
$$\aligned &e_1 = 28(z_1^6z_2^2 - z_1^2z_2^6),\\
           &e_3 = 56(z_1^7z_2 - z_1^5z_2^3 - z_1^3z_2^5 + z_1z_2^7),\\
           &e_5 = 8(z_1^8z_2 - 7z_1^5z_2^3 + 7z_1^3z_2^5 - z_1z_2^7),\\
           &e_7 = z_1^8 + z_2^8,\\
           &e_9 = 70z_1^4z_2^4,\\
           &a_1 = z_1^4 +z_2^4,\\
           &a_3 = z_1^4 - z_2^4,\\
           &a_5 = 4(z_1^3z_2 + z_1z_2^3).
  \endaligned
  \qquad
  \aligned &e_2 = 56(z_1^7z_2 + z_1^5z_2^3 - z_1^3z_2^5 - z_1z_2^7),\\
           &e_4 = z_1^8 - z_2^8,\\
           &e_6 = 8(z_1^7z_2 + 7z_1^5z_2^3 + 7z_1^3z_2^5 + z_1z_2^7),\\
           &e_8 = 28(z_1^6z_2^2 + z_1^2z_2^6),\\
           &a_0 = 1,\\
           &a_2 = 6z_1^2z_2^2,\\
           &a_4 = 4(z_1^3z_2 - z_1z_2^3),\\
           &
  \endaligned
$$
Let $(x,s) = (x_1, \dots , x_9, s_0, s_1, \dots , s_5)$ be
the corresponding coordinates in
$V(8) \dotplus V(0) \dotplus V(4)$.
\par
We have
$$\split \delta_{\lambda}(x,s) = &Q_1(x,s)(z_1^4 +z_2^4)
  + Q_2(x,s)6z_1^2z_2^2 + Q_3(x,s)(z_1^4 -z_2^4) \\
  &+ Q_4(x,s)4(z_1^3z_2 - z_1z_2^3) + Q_5(x,s)4(z_1^3z_2 +z_1z_2^3).
  \endsplit
$$
Direct calculations give us:
$$ \aligned
  Q_1(x,s) =
      &q_1(x) + x_7s_1 + x_9s_1
        + 6x_8s_2 - x_4s_3 + 8x_5s_4 \\
      & + 24x_2s_4 - 8x_6s_5 - 24x_3s_5
        + \epsilon(6s_1s_2 - 12s_4^2 - 12s_5^2) + s_0s_1, \\
  Q_2(x,s) =
      &q_2(x) + 2x_8s_1 + 6x_9s_2
        - 2x_1s_3 - 8x_5s_4 - 8x_2s_4 \\
      & - 8x_6s_5 + 8x_3s_5
        + \epsilon(2s_1^2 - 6s_2^2 - 2s_3^2 - 4s_4^2 + 4s_5^2) + s_0s_2, \\
  Q_3(x,s) =
      &q_3(x) + x_4s_1 + 6x_1s_2 - x_7s_3 + x_9s_3 \\
      & + 32x_3s_4 - 32x_2s_5
        + \epsilon(6s_2s_3 - 12s_4s_5) + s_0s_3, \\
  Q_4(x,s) =
      &q_4(x) + 2x_5s_1 + 6x_2s_1 - 6x_5s_2 + 6x_2s_2 - 8x_3s_3 \\
      & + 4x_8s_4 - 4x_9s_4 - 4x_1s_5
        +\epsilon(-3s_1s_4 - 3s_2s_4 + 3s_3s_5) + s_0s_4, \\
  Q_5(x,s) =
      &q_5(x) + 2x_6s_1 + 6x_3s_1 + 6x_2s_2 - 6x_3s_2 -8x_2s_3 \\
      & + 4x_1s_4 - 4x_8s_5 - 4x_9s_5
        + \epsilon(3s_1s_5 - 3s_2s_5 - 3s_3s_4) + s_0s_5,
   \endaligned \tag{1.2}
$$
where
$$ \aligned
    q_1(x) =
      &6x_7x_8 + 90x_8x_9 - 6x_4x_1 - 192x_5^2 - 96x_5x_2 - 192x_6^2 \\
      & - 96x_6x_3 + 384x_2^2 + 384x_3^2, \\
    q_2(x) =
      &2x_7^2 - 16x_8^2 - 50x_9^2 - 2x_4^2 - 64x_5^2 + 96x_5x_2 + 64x_6^2 \\
      & - 96x_6x_3 + 16x_1^2 + 128x_2^2 - 128x_3^2, \\
    q_3(x) =
      &-6x_7x_1 + 6x_8x_4 + 90x_9x_1 +48x_5x_6 - 336x_5x_3 \\
      & - 336x_6x_2 + 624x_2x_3, \\
    q_4(x) =
      &-3x_7x_5 - 21x_7x_2 + 12x_8x_5 - 132x_8x_2 + 15x_9x_5 \\
      & - 15x_9x_2 + 3x_4x_6 + 21x_4x_3 + 42x_6x_1 + 78x_1x_3, \\
    q_5(x) =
      &3x_7x_6 + 21x_7x_3 + 12x_8x_6 - 132x_8x_3 - 15x_9x_6 \\
      & + 15x_9x_3 - 3x_4x_5 - 21x_4x_2 + 42x_5x_1 + 78x_1x_2.
   \endaligned
$$

\par
\S \bf 2. \rm In this section we recall some facts about
$(G,G')$-sections.
\par
Let $G$ be a linear algebraic group, $X$ be an irreducible
quasiprojective variety, $G:X$ be a regular action, and
$G' \subset G$ be a subgroup of $G$.
\proclaim{Definition 2.1} An irreducible subvariety $X' \subset X$
is called $(G,G')$-section of $X$ iff
\roster
\item $\overline{G \cdot X'} = X$,
\item $G' \cdot X' = X'$,
\item $(G \cdot x') \cap X' = G' \cdot x'$ for all $x' \in X'$.
\endroster
\endproclaim
\par
Suppose $X'$ is $(G,G')$-section of $X$; then we have
the following isomorphism of the fields:
$$\C(X)^G \approx \C(X')^{G'}, \ \ \ \
  f \mapsto f \mid_{X'}.
$$
\par
Let $X'$ be $(G,G')$-section of $X$, $Y$ be an irreducible
quasiprojective variety, $G:Y$ be a regular action,
$F:Y \rightarrow X$ be a dominant $G$-morphism, and
$Y' \subset Y$ be an irreducible component of $F^{-1}(X')$.
\proclaim{Proposition 2.2} Suppose that $G' \cdot Y' = Y'$
and $F(Y')$ is dense in $X'$; then $Y'$ is
$(G,G')$-section of $Y$.
\endproclaim
\demo{Example 2.3} Let $G$ be a reductive linear algebraic
group, $G:X$ be a linear representation, and $H \subset G$ be
a stationary subgroup of general position
of the representation $G:X$. There exists an open
nonempty $G$-invariant subset $X^0$ such that $G_x$ conjugate
to $H$ for all $x \in X^0$. We have:
$$(X^H)^0 = (X^H) \cap X^0 = \{ x \in X^H \mid G_x = H \}
$$
is $(G,N(H))$-section of $X$, where $N(H)$ is the normalizer
of the subgroup $H$ in $G$.
\enddemo
\demo{Example 2.4} Consider the linear representation $PSL_2 : V(4)$.
As is known, the stationary subgroup of general
position of this representation is
$H = \{ e, \omega, \rho, \omega\rho \}$, where
$$e = \overline{\left( \matrix
         1  &  0  \\
         0  &  1  \endmatrix \right)} , \quad
  \omega = \overline{\left( \matrix
         0  &  1  \\
        -1  &  0  \endmatrix \right)} , \quad
  \rho = \overline{\left( \matrix
        -i  &  0  \\
         0  &  i  \endmatrix \right)} .
$$
It can easily be checked that
$N(H) = \langle \tau, \sigma \rangle$, where
$$\tau = \overline{\left( \matrix
        \theta^{-1}  &      0       \\
            0        &   \theta     \endmatrix \right)} , \quad
  \sigma = \overline{\frac {1}{\sqrt{2}}
  \left( \matrix
        \theta^3     &   \theta^7   \\
        \theta^5     &   \theta^5   \endmatrix \right)} , \quad
  \theta = exp(2 \pi i / 8) .
$$
We have $N(H) \simeq S_4$ and $N(H)/H \simeq S_3$.
It follows from Example 2.3 that
$$(V(4)^H)^0 = \{ f \in V(4)^H \mid (PSL_2)_f = H \}
$$
is $(PSL_2, N(H))$-section of $V(4)$.
\enddemo

\par
\S \bf 3. \rm In this section we construct
$(PSL_2, N(H))$-section $PX_{\lambda}^0$ of the
variety $PU_{\lambda}$ (see the definition of
$N(H)$ in \S 2).
\par
For the sake of convenience we write down how the groups
$H$ and $N(H)$ act in the space
$V(8) \dotplus V(0) \dotplus V(4)$:
$$\aligned \omega \cdot (x,s) =
      &(-x_1,x_2,-x_3,-x_4,x_5,-x_6,x_7,x_8,x_9,
        s_0,s_1,s_2,-s_3,s_4,-s_5), \\
           \rho \cdot (x,s) =
      &(x_1,-x_2,-x_3,x_4,-x_5,-x_6,x_7,x_8,x_9,
        s_0,s_1,s_2,s_3,-s_4,-s_5), \\
           \tau \cdot (x,s) =
      &(-x_1,-ix_3,-ix_2,x_4,-ix_6,-ix_5,x_7,-x_8,x_9, \\
      &  s_0,-s_1,s_2,-s_3,is_5,is_4), \\
           \sigma \cdot (x,s) =
      &(4x_3,-\frac{i}{4}x_1,ix_2,-8x_6,-\frac{i}{8}x_4,-ix_5, \\
      &\frac{1}{8}x_7 + \frac{7}{2}x_8 + \frac{35}{8}x_9,
       -\frac{1}{8}x_7 - \frac{1}{2}x_8 + \frac{5}{8}x_9,
       \frac{1}{8}x_7 - \frac{1}{2}x_8 + \frac{3}{8}x_9, \\
      &s_0,-\frac{1}{2}s_1 - \frac{3}{2}s_2,
       \frac{1}{2}s_1 - \frac{1}{2}s_2,2s_5,\frac{i}{2}s_3,-is_4).
  \endaligned \tag{3.1}
$$
We have
$$\aligned &(V(8) \dotplus V(0) \dotplus V(4))^H =
      \langle e_7, e_8, e_9, a_0, a_1, a_2 \rangle \\
           &(V(8) \dotplus V(0) \dotplus V(4))^{N(H)} =
      \langle 5e_7 + e_9, a_0 \rangle .
  \endaligned
$$
The decomposition of $N(H)$-module
$V(8) \dotplus V(0) \dotplus V(4)$ is as follows:
$$\aligned V(8) \dotplus V(0) \dotplus V(4) =
       &\langle e_1,e_2,e_3 \rangle \dotplus
        \langle e_4,e_5,e_6 \rangle \dotplus
        \langle e_8,7e_7 - e_9 \rangle \dotplus \\
       &\langle 5e_7 + e_9 \rangle \dotplus
        \langle a_0 \rangle \dotplus
        \langle a_1,a_2 \rangle \dotplus
        \langle a_3,a_4,a_5 \rangle .
  \endaligned
$$
\par
Set
$$p:V(8) \dotplus V(0) \dotplus V(4), \quad
  f_8 + f_0 + f_4 \mapsto  f_4.
$$
First we construct $(PSL_2, N(H))$-section $X_{\lambda}^0$ of
the variety $U_{\lambda}$ by applying Proposition 2.2
to $PSL_2$-morphism $p\mid_{U_{\lambda}}$ and
$(PSL_2, N(H))$-section $(V(4)^H)^0$ of $V(4)$
(see Example 2.4).
\proclaim{Lemma 3.1} $5e_7 + e_9 \in U_{\lambda}$.
\endproclaim
\demo{Proof} Consider the plane $\langle a_0, 5e_7 + e_9 \rangle$.
We have $N(H) \cdot \delta_{\lambda}(x,s) =
\delta_{\lambda}(N(H) \cdot (x,s)) =
\delta_{\lambda}(x,s)$ for all
$(x,s) \in \langle a_0, 5e_7 + e_9 \rangle$
(see (3.1)). Therefore,
$\delta_{\lambda}(\langle a_0, 5e_7 + e_9 \rangle)
 \subset V(4)^{N(H)} = \{ 0 \}$ and
$\langle a_0, 5e_7 + e_9 \rangle \subset
 \delta_{\lambda}^{-1}(0)$. Note also that
$a_0 \in U_{\lambda}$ and $a_0$ is a regular point
of $\delta_{\lambda}^{-1}(0)$. It follows that
$\langle a_0, 5e_7 + e_9 \rangle \subset U_{\lambda}$
and hence $5e_7 + e_0 \in U_{\lambda}$.
\enddemo
Consider
$\widetilde{X}_{\lambda} =
 p^{-1}(V(4)^H) \cap \delta_{\lambda}^{-1}(0)$.
{}From (1.2) and (3.1) we obtain the following equations
of the subvariety
$\widetilde{X}_{\lambda} \subset
 V(8) \dotplus V(0) \dotplus V(4)$:
$$\gathered
    s_3 = s_4 = s_5 = 0, \\
    q_1(x) + x_7s_1 + x_9s_1 + 6x_8s_2
      + \epsilon 6s_1s_2 + s_0s_1 = 0, \\
    q_2(x) + 2x_8s_1 + 6x_9s_2
      + \epsilon(2s_1^2 -6s_2^2) + s_0s_2 = 0, \\
    q_3(x) + x_4s_1 + 6x_1s_2 = 0, \\
    q_4(x) + 2x_5s_1 + 6x_2s_1 - 6x_5s_2 + 6x_2s_2 = 0, \\
    q_5(x) + 2x_6s_1 + 6x_3s_1 + 6x_6s_2 - 6x_3s_2 = 0.
  \endgathered \tag{3.2}
$$
\proclaim{Lemma 3.2} \roster
\item $5e_7 + e_9$ is a regular point of the subvariety
  $\widetilde{X}_{\lambda}$,
  $\dim T_{5e_7 + e_9}(\widetilde{X}_{\lambda}) = 7$.
\item Exactly one irreducible component of the subvariety
  $\widetilde{X}_{\lambda}$ (denote it by $X_{\lambda}$)
  contains $5e_7 + e_9$, $\dim X_{\lambda} = 7$.
\item $N(H) \cdot X_{\lambda} = X_{\lambda}$.
                     \endroster
\endproclaim
\demo{Proof} The proof of (1) is by direct calculations.
\par (2) is the consequence of (1).
\par (3). Since
$N(H) \cdot \widetilde{X}_{\lambda} = \widetilde{X}_{\lambda}$,
$N(H) \cdot (5e_7 + e_9) = 5e_7 + e_9$, and $5e_7 + e_9$
is a regular point of the subvariety
$\widetilde{X}_{\lambda}$, we see that
$N(H) \cdot X_{\lambda} = X_{\lambda}$.
\enddemo
It follows from Lemma 3.2 that $X_{\lambda}$ is an
irreducible component of the subvariety
$p^{-1}(V(4)^H) \cap U_{\lambda}$.
\par
We set
$$X_{\lambda}^0 = \{ (x,s) \in X_{\lambda} \mid
  p(x,s) \in (V(4)^H)^0 \}.
$$
Since $N(H) \cdot X_{\lambda} = X_{\lambda}$,
$N(H) \cdot (V(4)^H)^0 = (V(4)^H)^0$, we see that
$N(H) \cdot X_{\lambda}^0 = X_{\lambda}^0$.
It follows from Lemma 3.2 that $X_{\lambda}^0$
is a nonempty open subset of $X_{\lambda}$ and
$p(X_{\lambda}^0)$ is dense in $(V(4)^H)^0$.
This and Proposition 2.2 imply that
$X_{\lambda}^0$ is $(PSL_2, N(H))$-section of $U_{\lambda}$.
\par
Consider
$PX_{\lambda}^0 \subset PX_{\lambda} \subset PU_{\lambda}$.
It follows from the previous paragraph that $PX_{\lambda}^0$
is $(PSL_2, N(H))$-section of $PU_{\lambda}$.
We have the isomorphism of the fields
$$\C(PU_{\lambda})^{PSL_2} \approx
  \C(PX_{\lambda}^0)^{N(H)} \approx
  \C(PX_{\lambda})^{N(H)}.
$$
\par
Our goal now is to prove rationality of
$\C(PX_{\lambda})^{N(H)}$.
Note that $PX_{\lambda}$ is uniquely defined by the following
conditions:
\roster
\item $\overline{5e_7 + e_9} \in \p X_{\lambda}$,
\item $PX_{\lambda}$ is an irreducible component of
      $P\widetilde{X}_{\lambda}$,
\item the equations of the subvariety
      $P\widetilde{X}_{\lambda} \subset
       P(V(8) \dotplus V(0) \dotplus V(4))$
      are (3.2).
\endroster
(see Lemma 3.2)

\par
\S \bf 4. \rm In this section we define the linear space $R$,
the linear representation $N(H) : R$, the projective
representation $N(H) : \p ^8$, and $6$-dimensional irreducible
$N(H)$-invariant closed subvariety
$Y_{\lambda} \subset R \times \p ^8$ such that
$\C(PX_{\lambda})^{N(H)} \approx \C(Y_{\lambda})^{N(H)}$ and $H$ acts
on $Y_{\lambda}$ trivially.
\par
Let $R$ be $3$-dimensional linear space, and let
$r = (r_1,r_2,r_3)$ be coordinates in $R$. Define
the linear representation $N(H) : R$ in the following
way:
$$\tau \cdot (r_1,r_2,r_3) = (-r_1,r_3,r_2), \quad
  \sigma \cdot (r_1,r_2,r_3) = (-2r_3,r_1/2,-r_2).
$$
The subgroup $H \subset N(H)$ acts on $R$ trivially.
\par
Let $\overline{y} = (y_1:y_2:y_3:y_7:y_8: \dots :y_{12})$
be the homogeneous coordinates in $\p ^8$.
Define the projective representation $N(H) : \p ^8$
in the following way:
$$\aligned \tau \cdot \overline{y} =
     (&y_1:-y_3:-y_2:y_7:-y_8:y_9:y_{10}:-y_{11}:y_{12}), \\
           \sigma \cdot \overline{y} =
     (&\frac{1}{16}y_3:-16y_1:-y_2:
        \frac{1}{8}y_7 + \frac{7}{2}y_8 + \frac{35}{8}y_9:
        -\frac{1}{8}y_7 - \frac{1}{2}y_8 + \frac{5}{8}y_9: \\
      &\frac{1}{8}y_7 - \frac{1}{2}y_8 + \frac{3}{8}y_9:
        y_{10}:-\frac{1}{2}y_{11} - \frac{3}{2}y_{12}:
        \frac{1}{2}y_{11} - \frac{1}{2}y_{12}).
  \endaligned
$$
The subgroup $H \subset N(H)$ acts on $\p ^8$ trivially.
We have the regular action $N(H) : R \times \p ^8$.
The subgroup $H \subset N(H)$ acts on $R \times \p ^8$
trivially. Define the open $N(H)$-invariany subset
$${\p ^8}' = \{ \overline{y} \in \p ^8 \mid y_1y_2y_3 \neq 0 \}.
$$
\par
Set
$$M' = \{ (x,s) \in V(8) \dotplus V(0) \dotplus V(4) \mid
  s_3 = s_4 = s_5 = 0,\ x_1x_2x_3 \neq 0 \}.
$$
We see that $N(H) \cdot M' = M'$, $M = \overline{M'}$
is a linear subspace of
$V(8) \dotplus V(0) \dotplus V(4)$.
\par
Define the morphism
$$\gathered \pi : PM' \rightarrow R \times {\p ^8}', \\
      (x,s) \mapsto ((\frac{x_4}{x_1},\frac{x_5}{x_2},\frac{x_6}{x_3}),
      (\frac{x_2x_3}{x_1}:\frac{x_3x_1}{x_2}:\frac{x_1x_2}{x_3}):
       x_7:x_8:x_9:s_0:s_1:s_2)).
  \endgathered
$$
It can easily be checked that $\pi$ is $N(H)$-morphism and
fibers of $\pi$ are $H$-orbits.
\par
Note that
$PX_{\lambda} \subset P\widetilde{X}_{\lambda} \subset PM$.
Set
$$X'_{\lambda} = X_{\lambda} \cap M', \quad
  \widetilde{X}'_{\lambda} = \widetilde{X}_{\lambda} \cap M'.
$$
\proclaim{Lemma 4.1} $X'_{\lambda} \neq \emptyset$.
\endproclaim
Set
$$x^0 = 13i(5e_7 + e_9) + 5(4e_1 - ie_2 +e_3).
$$
Lemma 4.1 is a corollary of the following fact.
\proclaim{Lemma 4.2} $x^0 \in X_{\lambda}$.
\endproclaim
\demo{Proof} Consider the subgroup
$(\sigma) = \{ \sigma, \sigma^2, \sigma^3 = 1 \} \subset N(H)$.
We have
$$\aligned
     &V(8)^{(\sigma)} = \langle 5e_7 + e_9, 8e_4 -ie_5 - e_6,
          4e_1 - ie_2 + e_3 \rangle , \\
     &V(4)^{(\sigma)} = \langle 2(z_1^4 - z_2^4)
          + 4(z_1^3z_2 + z_1z_2^3)
          + 4i(z_1^3z_2 - z_1z_2^3) \rangle.
  \endaligned
$$
It follows from above that
$$\aligned \delta_{\lambda}(\alpha_1
     &(5e_7 + e_9)
       + \alpha_2(8e_4 - ie_5 - e_6)
       + \alpha_3(4e_1 - ie_2 + e_3)) \\
     &=q(\alpha_1,\alpha_2,\alpha_3)(2(z_1^4 - z_2^4)
       + 4(z_1^3z_2 + z_1z_2^3) + 4i(z_1^3z_2 -z_1z_2^3)).
  \endaligned \tag{4.1}
$$
Direct calculations give us
$$q(\alpha_1,\alpha_2,\alpha_3)
    = 24(5\alpha_1 \alpha_3 + i\alpha_2^2 - 13i\alpha_3^2). \tag{4.2}
$$
Consider $V(8)^{(\sigma)} \cap \widetilde{X}_{\lambda}$.
{}From (4.1) and (4.2) it follows that
$x^0, 5e_7 + e_9 \in V(8)^{(\sigma)} \cap \widetilde{X}_{\lambda}$
and $V(8)^{(\sigma)} \cap \widetilde{X}_{\lambda}$
is irreducible. On the other hand $5e_7 + e_9$ is a regular
point of $X_{\lambda}$ (Lemma 3.2). Hence
$V(8)^{(\sigma)} \cap \widetilde{X}_{\lambda} \subset X_{\lambda}$
and so $x^0 \in X_{\lambda}$.
\enddemo
{}From Lemma 4.1 it follows that $X'_{\lambda}$ is an open
nonempty $N(H)$-invariant subset of $X_{\lambda}$. We have
the isomorphism of the fields
$$\C(PX_{\lambda})^{N(H)}
  \approx \C(PX'_{\lambda})^{N(H)}. \tag{4.3}
$$
Notice that $PX'_{\lambda}$ is an irreducible component
of $P\widetilde{X}'_{\lambda}$ and
$\overline{x^0} \in PX_{\lambda}$.
\par
We have the isomorphism of the fields
$$\C(PX'_{\lambda})^{N(H)}
  \approx \C(\pi(PX'_{\lambda}))^{N(H)}. \tag{4.4}
$$
Notice that $\pi(PX'_{\lambda})$ is an irreducible component
of $\pi(P\widetilde{X}'_{\lambda})$, and
$$\pi(x^0) = ((0,0,0),(-5/4:20:-20:65:0:13:0:0:0))
  \in \pi(PX'_{\lambda}).
$$
It is not hard to obtain from (3.2) that
the equations of the subvariety
$\pi(P\widetilde{X}'_{\lambda}) \subset R \times {\p ^8}'$
are
$$\aligned 6y_7y_8
       &+ 90y_8y_9 +(-192r_3^2 - 96r_3 + 384)y_1y_2
        + (-192r_2^2 - 96r_2 + 384)y_1y_3 \\
       &+ (-6r_1)y_2y_3 + y_7y_{10} + y_9y_{11} + 6y_8y_{12}
        +6\epsilon y_{11}y_{12} + y_{10}y_{11} = 0, \\
           2y_7^2 -
       &16y_8^2 - 50y_9^2 + (64r_3^2 - 96r_3 - 128)y_1y_2
        +(-64r_2^2 + 96r_2 +128)y_1y_3 \\
       &+(-2r_1^2 + 16)y_2y_3 + 2y_8y_{11} + 6y_9y_{12}
        +\epsilon(2y_{11}^2 - 6y_{12}^2) + y_{10}y_{12} = 0 \\
           (48r_2r_3
       &- 336r_2 - 336r_3 + 624)y_1 -6y_7 +6r_1y_8 + 90y_9
        + r_1y_{11} + 6y_{12} = 0, \\
           (3r_1r_3
       &+ 21r_1 + 42r_3 + 78)y_2
        + (-3r_2 - 21)y_7 + (12r_2 - 132)y_8 \\
       &+ (15r_2 - 15)y_9 + (2r_2 + 6)y_{11} + (-6r_2 + 6)y_{12} = 0, \\
           (-3r_1
       &r_2 - 21r_1 + 42r_2 + 78)y_3 + (3r_1 + 21)y_7 + (12r_3 -132)y_8 \\
       &+(-15r_3 + 15)y_9 + (2r_3 + 6)y_{11} + (6r_3 - 6)y_{12} = 0.
  \endaligned \tag{4.5}
$$
\par
Equations (4.5) define in $R \times \p ^8$ the subvariety
$\widetilde{Y}_{\lambda}$. The closure of
$\pi(P\widetilde{X}'_{\lambda})$ in $R \times \p ^8$
is a union of some irreducible components of
$\widetilde{Y}_{\lambda}$. Let $Y_{\lambda}$ be a closure of
$\pi(PX'_{\lambda})$ in $R \times \p ^8$.
We see that
$N(H) \cdot \widetilde{Y}_{\lambda} = \widetilde{Y}_{\lambda},
\ N(H) \cdot Y_{\lambda} = Y_{\lambda}$,
and $Y_{\lambda}$ is
an irreducible component of the subvariety
$\widetilde{Y}_{\lambda}$. We have the isomorphism of
the fields
$$\C(\pi(PX'_{\lambda}))^{N(H)} \approx \C(Y_{\lambda})^{N(H)}. \tag{4.6}
$$
{}From (4.3), (4.4), and (4.6) we obtain the isomorphism of
the fields
$$\C(PX_{\lambda})^{N(H)} \approx \C(Y_{\lambda})^{N(H)}.
$$
\par
Our goal now is to prove rationality of
$\C(Y_{\lambda})^{N(H)}$.
Note that the following
conditions hold for $Y_{\lambda}$:
\roster
\item $\pi(\overline{x^0}) \in Y_{\lambda}$,
\item $Y_{\lambda}$ is an irreducible component of
  $\widetilde{Y}_{\lambda}$,
\item the equations of the subvariety
  $\widetilde{Y}_{\lambda} \subset R \times \p ^8$ are (4.5).
\endroster

\par
\S \bf 5. \rm In this section we proove
rationality of $\C(Y_{\lambda})^{N(H)}$.
\par
Set
$$\gathered
    \eta : \widetilde{Y}_{\lambda} \rightarrow R,
    \quad (r,\overline{y}) \mapsto r, \\
    \beta : \widetilde{Y}_{\lambda} \rightarrow \p ^8,
    \quad (r,\overline{y}) \mapsto \overline{y}.
  \endgathered
$$
We have $\eta(\pi(\overline{x^0})) = 0$.
It follows from (4.5) that $\beta(\eta^{-1}(r))$
is an intersection of 2 quadrics and
3 hyperplanes in $\p ^8$.
\proclaim{Lemma 5.1} $\eta^{-1}(0)$ is irreducible
and $3$-dimensional.
\endproclaim
\demo{Proof} The proof is by direct calculations.
\enddemo
Set
$$R' = \{ r \in R \mid \eta^{-1}(r) \ is \ irreducible
  \ and \ 3-dimensional \}.
$$
{}From Lemma 5.1 it follows that $R'$ is an open nonempty
$N(H)$-invariant subset of $R, \ 0 \in R'$, and
$\eta^{-1}(R')$ is an open nonempty $N(H)$-invariant
subset of $Y_{\lambda}$. We have the isomorphism of
the fields
$$\C(Y_{\lambda})^{N(H)} \approx \C(\eta^{-1}(R'))^{N(H)}.
$$
\par
Let us prove rationality of $\C(\eta^{-1}(R'))^{N(H)}$.
\par
Consider the bundle
$$\eta \mid_{\eta^{-1}(R')} : \eta^{-1}(R') \rightarrow R'.
$$
This bundle has $N(H)$-section
$$r \mapsto (r,u'(r)), \quad u'(r) = (0:0:0:0:0:0:1:0:0).
$$
\proclaim{Lemma 5.2} There exists an open nonempty
$N(H)$-invariant subset $R'' \subset R$ such that
\roster
\item $R'' \ni 0$,
\item the bundle
$$\eta \mid _{\eta^{-1}(R'')} : \eta^{-1}(R'') \rightarrow R''
$$
has $N(H)$-section
$$r \mapsto (r,u''(r)) = (r,u_1''(r): \dots : u_9''(r))
$$
such that $u''_7(r) = u_8''(r) = u_9''(r) = 0$ for $r \in R''$,
\item $u''(0) = (-5/4:20:-20:65:0:13:0:0:0)$.
\endroster
\endproclaim
See the proof in \S 6.
\par
By (4.5) and lemma 5.2 it follows that
$$\langle u'(r),u''(r) \rangle
  \subset \beta(\eta^{-1}(r)) \quad for\ r \in R' \cap R''.
$$
\par
Set
$$\aligned N = \{
    &\overline{y} \in \p ^8 \mid
      y_1 = y_2 = y_3 = y_7 + 7y_9 = y_{10} = 0 \}, \\
          N(r) = \langle
    &u'(r), u''(r), (1:0:0:0: \dots ), (0:1:0:0: \dots ), \\
    &(0:0:1:0: \dots ) \rangle \subset \p ^8, \quad r \in R' \cap R''.
  \endaligned
$$
We have $N(H) \cdot N = N,\ g \cdot N(r) = N(g \cdot r)$
for $g \in N(H)$.
\proclaim{Lemma 5.3} There exists a nonempty open
$N(H)$-invariant subset
$0 \in R''' \subset (R' \cap R'')$
such that
\roster
\item $\dim N(r) = 4$,
\item $N(r) \cap N = \emptyset$
\endroster
for $r \in R'''$.
\endproclaim
\demo{Proof} From Lemma 5.2 it follows that
$\dim N(0) = 4, \ N(0) \cap N = \emptyset$.
{}From this it follows the lemma.
\enddemo
For $r \in R'''$ let
$$\gamma_r : P^8 \rightarrow N
$$
be the projection of $\p ^8$ to $N$ from $N(r)$.
\proclaim{Lemma 5.4} There exists a nonempty open
$N(H)$-invariant subset $0 \in R'''' \subset R'''$
such that $\gamma_r(\beta(\eta^{-1}(r))) = N$
for $r \in R''''$.
\endproclaim
\demo{Proof} It can easily be checked that
$\gamma_0(\beta(\eta^{-1}(0))) = N$.
{}From this it follows the lemma.
\enddemo
We have the isomorphism of the fields
$$\C(\eta^{-1}(R'))^{N(H)} \approx
  \C(\eta^{-1}(R''''))^{N(H)}.
$$
\par
Let us prove rationality of
$\C(\eta^{-1}(R''''))^{N(H)}$.
\par
Recall the following fact.
\proclaim{Lemma 5.5} Let $X \subset \p ^n$ be an intersection
of $5$-dimensional linear subspace and two quadrics, and
let $M_1, M_2 \subset \p ^n$ be linear subspaces.
Suppose $X$ is irreducible,
$\dim X = 3,\ M_1 \cap M_2 = \emptyset,\
 \dim M_1 = n - 4,\ \dim M_2 = 3,\ M_1 \cap X$
contains a line, and $p_2(X) = M_2$,
where $p_2$ is the projection of $\p ^n$ to
$M_2$ from $M_1$; then $p_2 \mid _X$ is a birational
isomorphism of $X$ and $M_2$.
\endproclaim
{}From Lemmas 5.4 and 5.5 it follows that
$$\gamma_r \mid _{\beta(\eta^{-1}(r))} :
  \beta(\eta^{-1}(r)) \rightarrow N
$$
is a birational isomorphism for $r \in R''''$.
Therefore,
$$\Gamma : \eta^{-1}(R'''') \rightarrow R'''' \times N,
  \quad (r,\overline{y}) \mapsto (r,\gamma_r(\overline{y}))
$$
is a birational $N(H)$-isomorphism. The birational
isomorphism $\Gamma$ defines the isomorphism of the fields
$$\C(\eta^{-1}(R''''))^{N(H)} \approx
  \C(R'''' \times N)^{N(H)}.
$$
Rationality of the field
$$\C(R'''' \times N)^{N(H)} \approx \C(R \times N)^{N(H)}
$$
is a consequence of Noname lemma and Castelnuovo's
therem.

\par
\S \bf 6. \rm In this section we prove Lemma 5.2.
\par
Let $X_1 \subset P(V(8) \dotplus V(0) \dotplus V(4))$
be the projectivization of
$\overline{PSL_2 \cdot \langle z_1^8,z_1^7z_2,z_1^6z_2^2 \rangle}$
and let $X_2 \subset P(V(8) \dotplus V(0) \dotplus V(4))$
be the projectivization of
$\overline{PSL_2 \cdot \langle 5e_7 + e_9 \rangle}$.
It is obvious that $X_1$ and $X_2$ are irreducible,
$\dim X_1 = \dim X_2 = 3$, and $f \in X_1$ iff $f$ has
a root of multiplicity $\ge 6$ (as element of $V(8)$).
\par
It is clear that
$\delta_{\lambda}(\langle z_1^8,z_1^7z_2,z_1^6z_2^2 \rangle) = 0$,
the differential
$d(\delta_{\lambda} \mid _{V(8)}) \mid _{z_1^6z_2^2}$
is surjective. This implies that $X_1$ is an irreducible
component of
$P(\delta_{\lambda}^{-1}(0) \cap V(8))$.
Note also that
$$\deg X_1 = 18
$$
(see \cite{12}).
\par
Since $\delta_{\lambda}(5e_7 + e_9) = 0$, the differential
$d(\delta_{\lambda} \mid _{V(8)}) \mid _{5e_7 + e_9}$
is sujective, we see that $X_2$ is an irreducible component of
$P(\delta_{\lambda}^{-1}(0) \cap V(8))$.
Since the stabilizer of $\overline{5e_7 + e_9}$ in
$PSL_2$ coinsides with $N(H)$ and $5e_7 + e_9$ has
distinct roots, we have
$$\deg X_2 = \frac{8 \cdot 7 \cdot 6}{N(H)} = 14.
$$
\par
{}From the considerations above we obtain the following fact.
\proclaim{Lemma 6.1}
$P(\delta_{\lambda}^{-1}(0) \cap V(8)) = X_1 \cup X_2$.
\endproclaim
For $r \in R$ define
$$L(r) = \{ \overline{(x,s)} \mid
  x_4 = r_1x_1,\ x_5 = r_2x_2,\ x_6 = r_3x_3 \}.
$$
We shall describe $L(r) \cap X_1$ and $L(r) \cap X_2$.
\par
Set
$$\gathered L_0 = \{\overline{(x,s)} \mid
    x_1 = x_2 = x_3 = x_4 = x_5 = x_6 = 0 \}, \\
  L_1(r) = \{ \overline{(x,s)} \mid
    x_1 \neq 0,\ x_4 = r_1x_1,\ x_2 = x_3 = x_5 = x_6 = 0 \}, \\
  L_2(r) = \{ \overline{(x,s)} \mid
    x_2 \neq 0,\ x_5 = r_2x_2,\ x_1 = x_3 = x_4 = x_6 = 0 \}, \\
  L_3(r) = \{ \overline{(x,s)} \mid
    x_3 \neq 0,\ x_6 = r_3x_3,\ x_1 = x_2 = x_4 = x_5 = 0 \}, \\
  \widetilde{L}_1(r) = \{ \overline{(x,s)} \mid
    x_2x_3 \neq 0,\ x_5 = r_2x_2,\ x_6 = r_3x_3,\ x_1 = x_4 = 0 \}, \\
  \widetilde{L}_2(r) = \{ \overline{(x,s)} \mid
    x_1x_3 \neq 0,\ x_4 = r_1x_1,\ x_6 = r_3x_3,\ x_2 = x_5 = 0 \}, \\
  \widetilde{L}_3(r) = \{ \overline{(x,s)} \mid
    x_1x_2 \neq 0,\ x_4 = r_1x_1,\ x_5 = r_2x_2,\ x_3 = x_6 = 0 \}, \\
  L^0(r) = \{ \overline{(x,s)} \mid
    x_1x_2x_3 \neq 0,\ x_4 = r_1x_1,\ x_5 = r_2x_2,\ x_6 = r_3x_3 \}.
  \endgathered
$$
The linear subspace $L(r)$ is the disjoint union of the subsets
$L_0,\ L^0(r),\ L_i(r)$,
$\widetilde{L}_i(r),\ i = 1,2,3$.
We have
$$\gathered g \cdot L(r) = L(g \cdot r), \quad
      g \cdot L^0(r) = L^0(g \cdot r), \quad
      g \cdot L_0 = L_0, \\
    g \cdot L_j(r) = L_{\kappa(g)(j)}(g \cdot r),
  \endgathered
$$
for $g \in N(H),\ r \in R$, where
$$\gathered \kappa : N(H) \rightarrow S_3, \\
    \kappa(\tau) = \left( \matrix
      1 & 2 & 3 \\
      1 & 3 & 2           \endmatrix \right) , \quad
    \kappa(\sigma) = \left( \matrix
      1 & 2 & 3 \\
      2 & 3 & 1             \endmatrix \right)
  \endgathered
$$
is the homomorphism of the groups.
\proclaim{Lemma 6.2} There exist an open nonempty $N(H)$-invariant
subset $0 \in R'' \subset R$ such that
$L(r) \cap P(\delta_{\lambda}^{-1}(0) \cap V(8))$ is $32$
points of multiplicity $1$ for $r \in R''$ such that
\roster
\item $\widetilde{L}_j(r) \cap X_l = \emptyset,
  \ 1 \le j \le 3,\ 1 \le l \le 2$;
\item $L_0 \cap X_1 = \emptyset,\ |L_0 \cap X_2| = 4$;
\item $|L_j(r) \cap X_l| = 2,\ 1 \le j \le 3,\ 1 \le l \le 2$;
\item $|L^0(r) \cap X_1| = 12,\ |L^0(r) \cap X_2| = 4$.
\endroster
\endproclaim
\demo{Proof} Set
$$\aligned R^0 = \{
    &r \in R \mid 48r_2r_3 - 336r_2 - 336r_3 + 624 \neq 0, \\
    &3r_1r_3 + 21r_1 + 42r_3 + 78 \neq 0,
     -3r_1r_2 - 21r_1 + 42r_2 + 78 \neq 0 \}.
  \endaligned
$$
{}From (1.2) it follows that if $r \in R^0,\ 1 \le j \le 3$, then
$\widetilde{L}_j(r) \cap
 P(\delta_{\lambda}^{-1}(0) \cap V(8)) = \emptyset$.
\par
It is sufficient to prove that
\par a)
$|L^0(0) \cap P(\delta_{\lambda}^{-1}(0) \cap V(8))| = 16$,
\par b)
$L_0 \cap X_1 = \emptyset,\ |L_0 \cap X_2| = 4,\ |L_j(r) \cap X_l| = 2
  \ (1 \le j \le 3,\ 1 \le l \le 2)$ for $r \in R$.
\par
Relation a) can be proved by strightforward calculations.
\par
Let us prove b).
\par
Consider $\overline{f} \in (L_1(r) \cup L_0) \cap PV(8)$.
If $(a:b)$ is a root of $f$ of multiplicity $m$, then
$(a:-b)$ is a root of $f$ of multiplicity $m$.
It follows that if $(a:b)$ is a root of $f$ of multiplicity
$\ge 6$, then $(a:b) = (1:0)$ or $(a:b) = (0:1)$.
Suppose $\overline{f} \in L_0$; then neither $(1:0)$ nor $(0:1)$
is a root of $f$ of multiplicity $\ge 6$. Therefore,
$$L_0 \cap X_1 = \emptyset. \tag{6.1}
$$
Suppose $\overline{f} \in L_1(r) \cap X_1$. If $(1:0)$ is a root of
$f$ of multiplicity $\ge 6$, then
$\overline{f} = \overline{-e_1 - r_1e_4 + r_1e_7 + e_8}$.
If $(0:1)$ is a root of $f$ of multiplicity $\ge 6$, then
$\overline{f} = \overline{e_1 + r_1e_4 + r_1e_7 + e_8}$.
It follows that
$$|L_1(r) \cap X_1| = 2. \tag{6.2}
$$
\par
Direct calculations give us
$$L_0 \cap P(\delta_{\lambda}^{-1}(0) \cap V(8)) =
  \{ \overline{5e_7 \pm e_9}, \overline{15e_7 \pm 5e_8 - e_9} \}.
  \tag{6.3}
$$
Taking into account (6.1) and (6.3), we obtain
$$ |L_0 \cap X_2| = 4.
$$
Direct calculations give us
$$\aligned L_1(r) \cap
  &P(\delta_{\lambda}^{-1}(0) \cap V(8)) =
   \{ \overline{\pm (e_1 + r_1e_4) + r_1e_7 + e_8}, \\
  &\overline{\pm (ae_1 + r_1ae_4) + (90 - 5r_1^2)e_7 - 5r_1e_8 + 6e_9} \},
  \endaligned \tag{6.4}
$$
where $a^2 = 25(r_1^2 - 36)$.
Using (6.2) and (6.4), we get
$$|L_1(r) \cap X_2| = 2.
$$
\par
We have
$$\sigma \cdot L_1(r) = L_2(\sigma \cdot r), \quad
  \sigma \cdot L_2(r) = L_3(\sigma \cdot r), \quad
  \sigma \cdot L_3(r) = L_1(\sigma \cdot r).
$$
For $2 \le j \le 3,\ 1 \le l \le 2$, we have
$$|L_j(r) \cap X_l| =
  |(\sigma^{1-j} \cdot L_j(r)) \cap (\sigma^{1-j} \cdot X_l)| =
  |L_1(\sigma^{1-j} \cdot r) \cap X_l| = 2.
$$
\enddemo
\proclaim{Corollary} $L^0(r) \cap X_2$ is $H$-orbit for
$r \in R''$.
\endproclaim
\demo{Proof} It is clear that a stabilizer of any
$\overline{x} \in L^0(r)$ in the group $H$ is trivial.
Therefore, any $H$-invariant finite subset of $L^0(r)$
of $4$ points is $H$-orbit. Hence $L^0(r) \cap X_2$
is $H$-orbit.
\enddemo
\demo{Proof of Lemma 5.2} Set
$$(r,u''(r)) = \pi(X_2 \cap L^0(r)).
$$
{}From Lemma 6.2 and Corollary of Lemma 6.2 follow
statements (1) and (2) of Lemma 5.2.
\par
Let us prove statement (3) of Lemma 5.2.
\par
It can easily be checked that $x^0$ has not a root
of multiplicity $\ge 6$ (as an element of $V(8)$).
{}From Lemma 6.1 it follows that
$\overline{x^0} \in X_2$. We have
$$u''(0) = u''(\pi(\overline{x^0})) =
  \pi(\overline{x^0}) = (-5/4:20:-20:65:0:13:0:0:0).
$$
\enddemo

\enddocument